\documentclass[onecolumn,aps,preprint,showpacs,amsmath,amssymb]{revtex4}
\usepackage{graphicx}
\usepackage{latexsym}
\usepackage{amsmath}
\usepackage{dcolumn} 

\usepackage{bm}

\begin{document}
\newcommand{\eq}{\begin{equation}}                                                                         
\newcommand{\eqe}{\end{equation}}             
 
\title{Lissajous curves  with a finite sum of prime number frequencies}

%\author{Imre F. Barna}
%\address{Wigner Research Center for Physics, 
%\\ Konkoly-Thege Mikl\'os \'ut 29 - 33, 1121 Budapest, Hungary \\
%Email: barna.imre@wigner.hu
%} 
\author{Imre F. Barna$^{1}$ and L. M\'aty\'as$^{2}$}
\address{ $^1$ Wigner Research Center for Physics, 
\\ Konkoly-Thege Mikl\'os \'ut 29 - 33, 1121 Budapest, Hungary \\ 
Email: barna.imre@wigner.hu\\
$^2$  Department of Bioengineering, Faculty of Economics, Socio-Human Sciences
and Engineering, Sapientia Hungarian University of Transylvania
Libert\u{a}tii sq. 1, 530104 Miercurea Ciuc, Romania} 

\date{\today}

%\maketitle
%%%%%%%%%%%%%%%%%%%%%%%%%%%%%%%%%%%%%%%%%%%%%%%%%%%%%%%%%%%%%%%%%%%%%%%
\begin{abstract} 
The Ulam spiral inspired us to calculate and present Lissajous curves where the orthogonally added functions are a finite sum of sinus and cosines 
functions with consecutive prime number frequencies. 
\end{abstract}

%\draft
%\pacs{47.10.A−, 47.10.ab, 47.55.Hd }
% mathemtical formulation,  conservation laws and constitutive relations, stratified flows, 
\maketitle,

%%%%%%%%%%%%%%%%%%%%%%%%%%%%%%%%%%%%%%%%%%
We may say that prime numbers fascinates mankind more than two thousand years. The scientific literature of 
number theory - which in great part deals with prime - is enormous and fills libraries. 
Number theory is not our field of interest at all, so it is not our duty to give any kind of 
overview of the field, therefore we just mention two works about primes \cite{schleich,hodge}. 
(Our scientific interest is laser-matter interaction \cite{laser} and 
self-similar solutions of non-linear partial differential 
equations of flow systems \cite{flow}.)   
We just would like to show a tiny idea about primes which might be interesting to the experts. 
There are two starting points which gave as the idea. 
The first is the Ulam spiral which was found by Stanislaw Ulam in 1967 \cite{ulamsp}. 
 The left figure  in Fig. (1) basically shows how the spiral is defined, the middle figure 
shows how the prime numbers are distributed among the first 400 natural numbers, and the last 
right figure presents the prime distribution on a much lager scales. 
The primes are represented with dark dots, spot and with short lines. 
It is evident, that there is a non-trivial correlation of primes even on large scales in this representation. 

Our second starting point is the definition of Lissajous (or Bowditch) curves  \cite{lawr} which are a well-known object for 
physicists. The parametric formula of the curve reads 
\begin{eqnarray}
x(t) &=& sin(a\cdot t + \delta),   \nonumber  \\  
y(t) &=& cos(b\cdot t), 
\end{eqnarray}
where 'a', 'b' are the relative frequencies and $ \delta$ is the phase between the two 
oscillations. 
Figure 2 presents three classical Lissajous curves with various relative frequencies and phases between the two trigonometric functions. 
 The length of the corresponding parametric curve is defined as 
\eq
L = \int_0^{2\pi} \sqrt{\dot{x}(t)^2 + \dot{y}(t)^2  } dt,  
\eqe
where prime means derivation with respect to the parameter $t$. 
With this formula it is trivial to get back the circumference of a unit circle. 
It is also clear that the length of various Lissajous curves are proportional to $\pi$. 
Differential geometrical analysis helps to derive and study additional parameters of the curve.
%%%%%%%%%%%%%%%%%%%%%%%%%%%%%%%%%%%%%%%%%%%%
\begin{figure} [!h]
 \hspace*{-0.6cm} 
\scalebox{0.7}{
\rotatebox{0}{\includegraphics{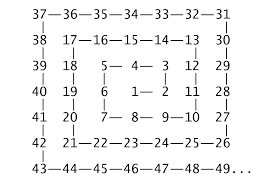}} } 
 \hspace*{-0.9cm} 
\scalebox{0.35}{
\rotatebox{0}{\includegraphics{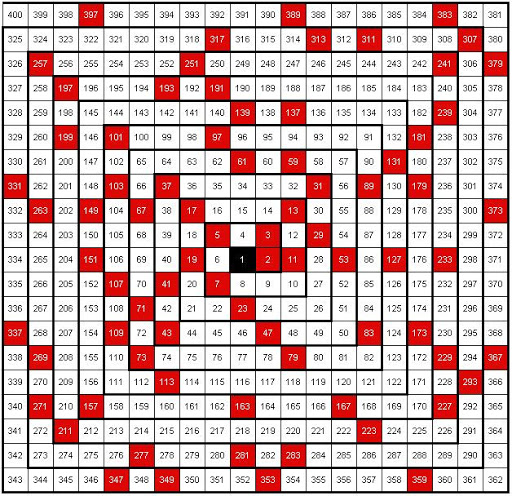}}}
 \scalebox{0.5}{
\rotatebox{0}{\includegraphics{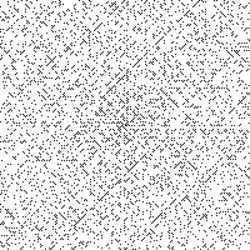}}}
\caption{The Ulam spiral. The left figure is just the definition of the spiral, the middle 
figure shows the prime distributions among the first 400 natural numbers, the right figure 
presents the large scale distribution of primes in the Ulam spiral.}	
\label{Ulam}       % Give a unique label
\end{figure}  
%%%%%%%%%%%%%%%%%%%%%%%%%%%%%%%%%%%%%%%%%%%%
\begin{figure} [!h]
\scalebox{0.75}{
\rotatebox{0}{\includegraphics{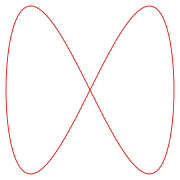}} } 
\scalebox{0.75}{
\rotatebox{0}{\includegraphics{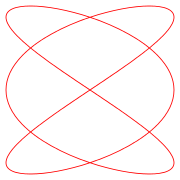}}}
\scalebox{0.75}{
\rotatebox{0}{\includegraphics{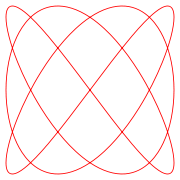}}}
\caption{Three classical Lissajous curves. The three parameter sets (from left to right) are $(a = 1,  \delta = \frac{\pi}{2}, b = 2) $, 
$ (a = 3,   \delta = \frac{\pi}{2},   b = 2),$ and  $(a = 3,  \delta = \frac{\pi}{4}, b =4)$  }	
\label{Liss}       % Give a unique label
\end{figure}
%%%%%%%%%%%%%%%%%%%%%%%%%%
%%%%%%%%%%%%%%%%%%%%%%%%%%%%%%%%%%%%%%%%%%%%
\begin{figure} [!h]
\scalebox{0.45}{
\rotatebox{0}{\includegraphics{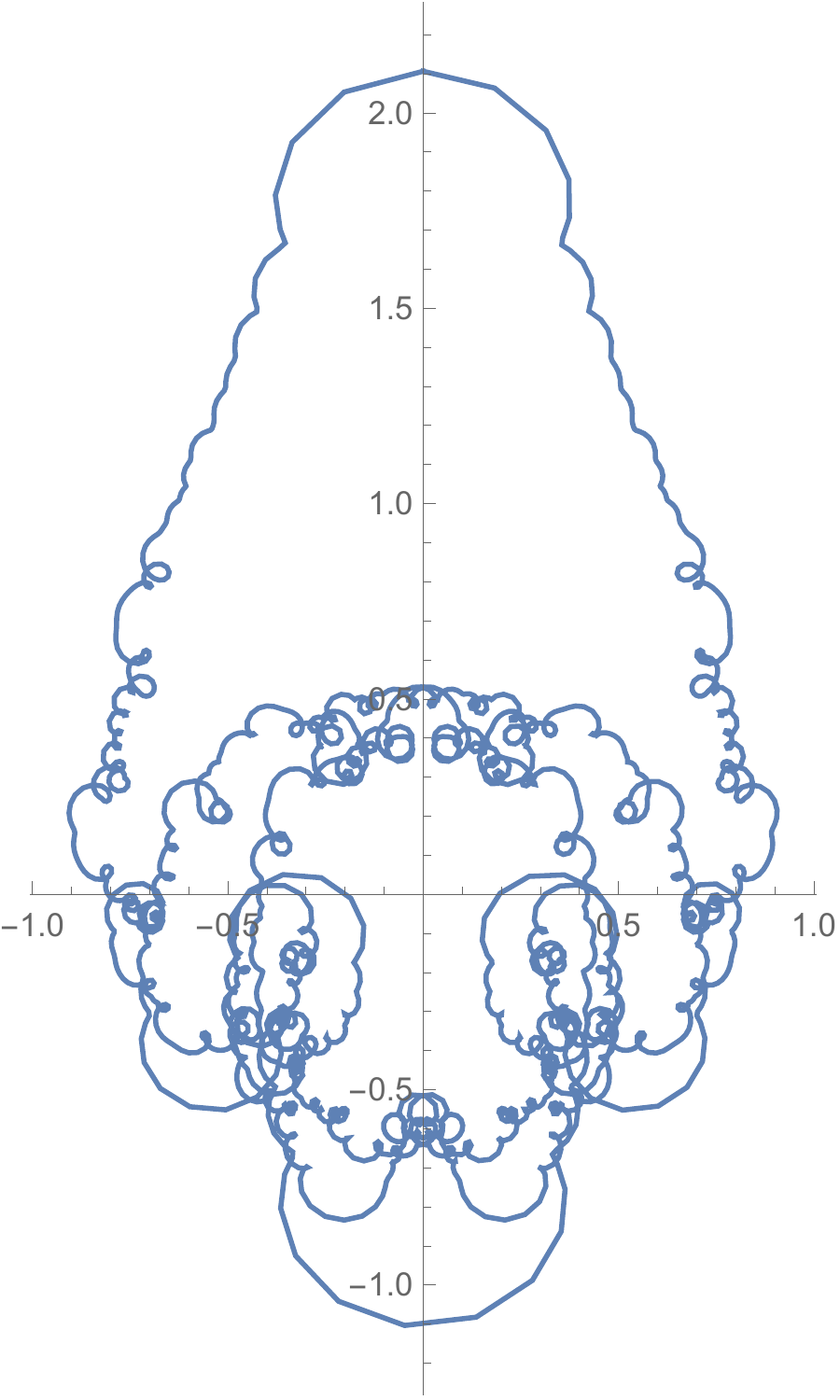}} } 
\scalebox{0.45}{
\rotatebox{0}{\includegraphics{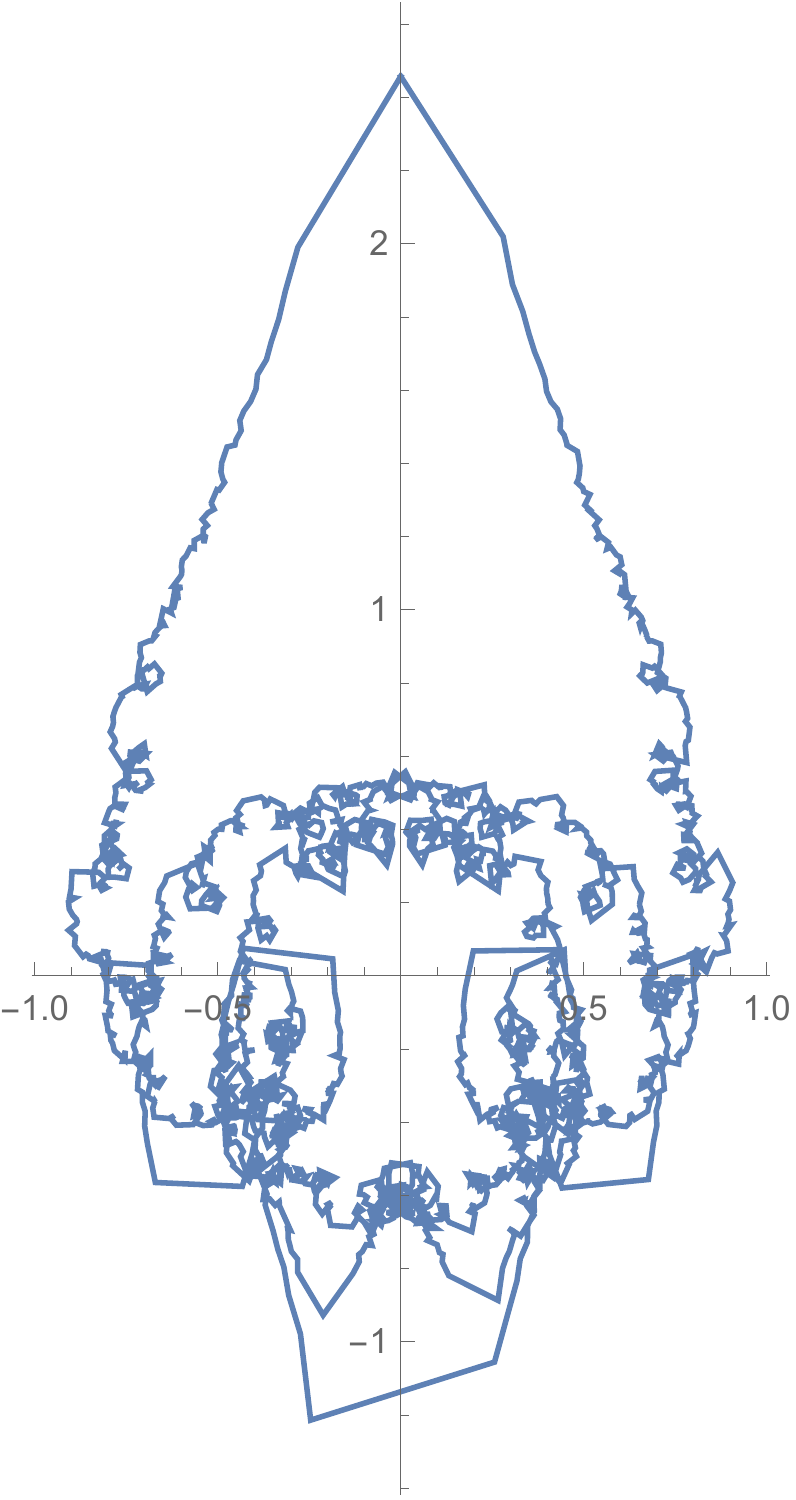}}}
\scalebox{0.45}{
\rotatebox{0}{\includegraphics{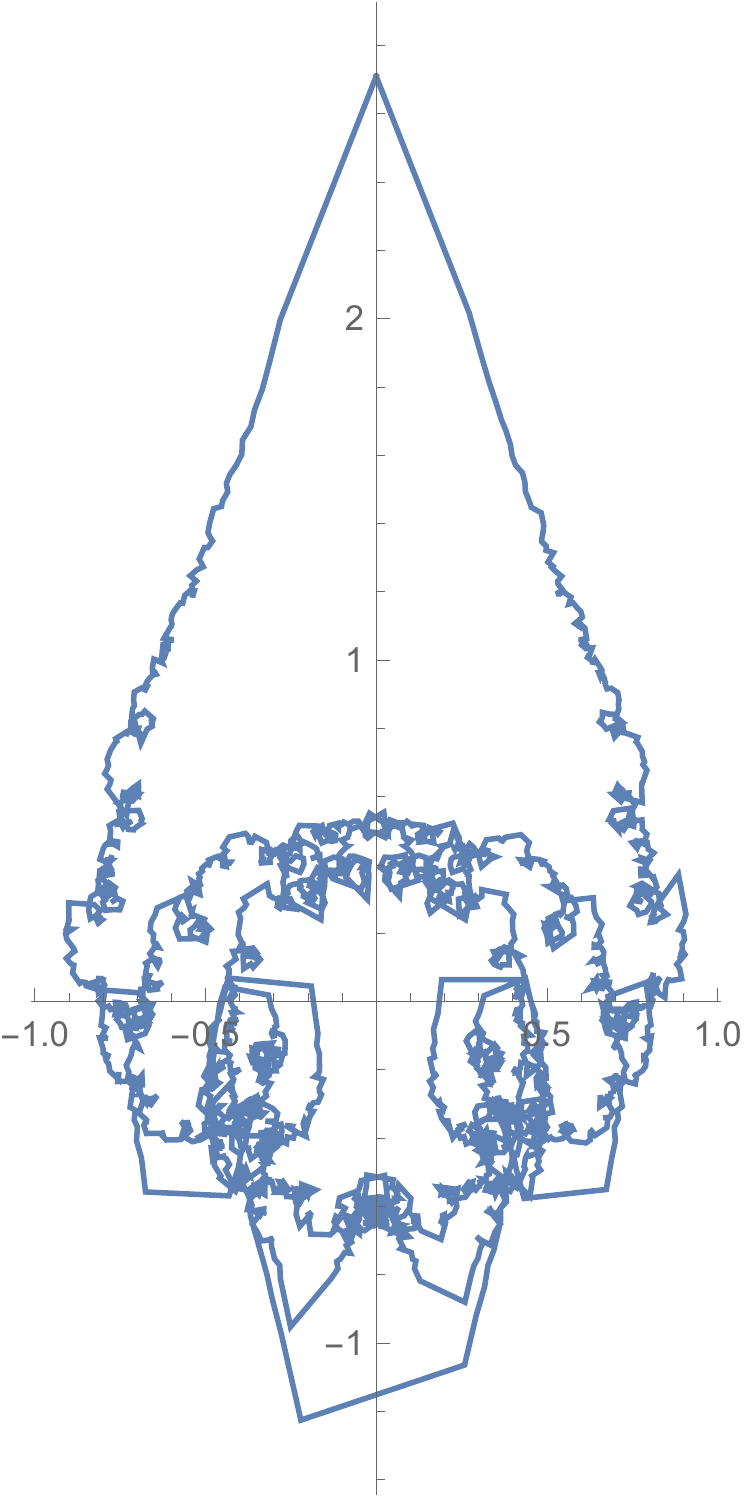}}}
\caption{Our Lissajous curves with different kind of finite Fourier sums with prime number frequencies. 
From left to right, the sum of the first 100, 1000 and 5000 primes were taken.}	
\label{Liss2}       % Give a unique label
\end{figure}
%%%%%%%%%%%%%%%%%%%%%% 
Let's try to image the distribution of primes somehow with the help  of the Lissajous curves. 

We applied the next parametric formula for the curve  
\begin{eqnarray}
x(t) &=&  \sum_{i=1}^N  \frac{ sin(a_i\cdot t )}{a_i},   \nonumber  \\  
y(t) &=& \sum_{i=1}^N \frac{cos(a_i\cdot t)}{a_i}, 
\end{eqnarray}
where $a_i$s are the first $N$ prime numbers. 
Figure 3 shows the Lissajous curves for $N = 100, 1000$ and $5000$ where the last  
primes are $541,7919$ and $104729$. 
To achieve a finite surface for the Lissajous curve we divide the sum of sinus and cosines  functions with 
the corresponding prime number.  Note, the slight left-right symmetry breaking of the curves. 
These are our starting points, of course the applied curves can be changed in numerous ways.
   
As second case, Figure 4 presents two curves where only the second neighbor primes are considered to the 'x' 
and 'y' coordinates  such as ${2,5,11}$ and ${3,7,13}$. 
Note, the much quicker convergence, it is not possible to see the differences between the two 
figures with naked eyes. We tried to modify Eq. (3) with additional logarithmic, square root or power law functions 
of the argument of the sinus and cosines function to create much more internal structure of the curves 
at larger number of primes. Unfortunately in vain. This is the present endpoint of our idea and analysis.        
(The presented calculations and figures were evaluated with Maple 12.)
It can happen that our toy model might give idea for such kind of further investigations.      
%%%%%%%%%%%%%%%%%%%%%%%%%%%%%%%%%%%%%%%%%%%%%%%%%%%%%%%%%%                                    
 \begin{figure} [!h]
\scalebox{0.45}{
\rotatebox{0}{\includegraphics{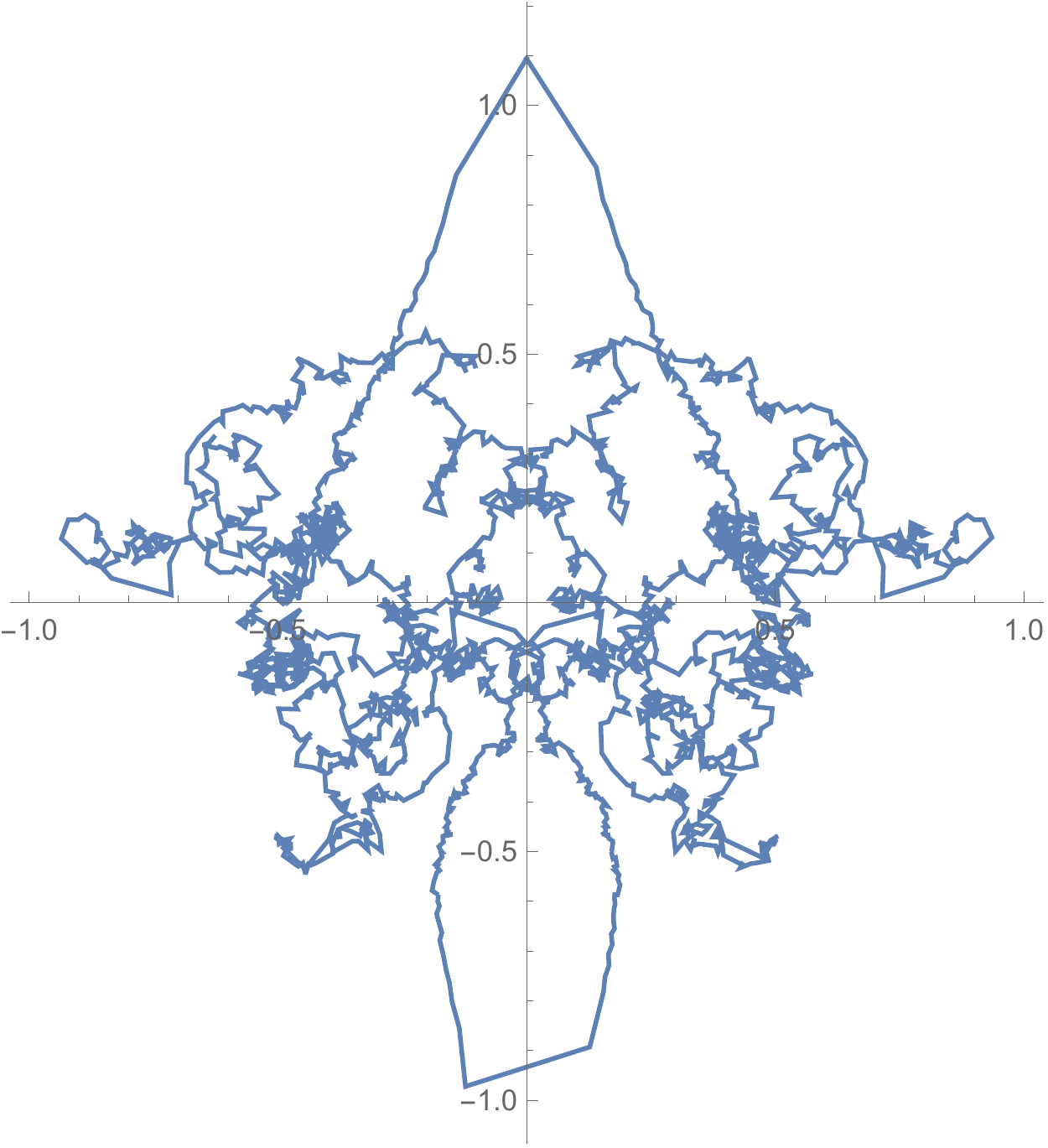}} } 
\scalebox{0.45}{
\rotatebox{0}{\includegraphics{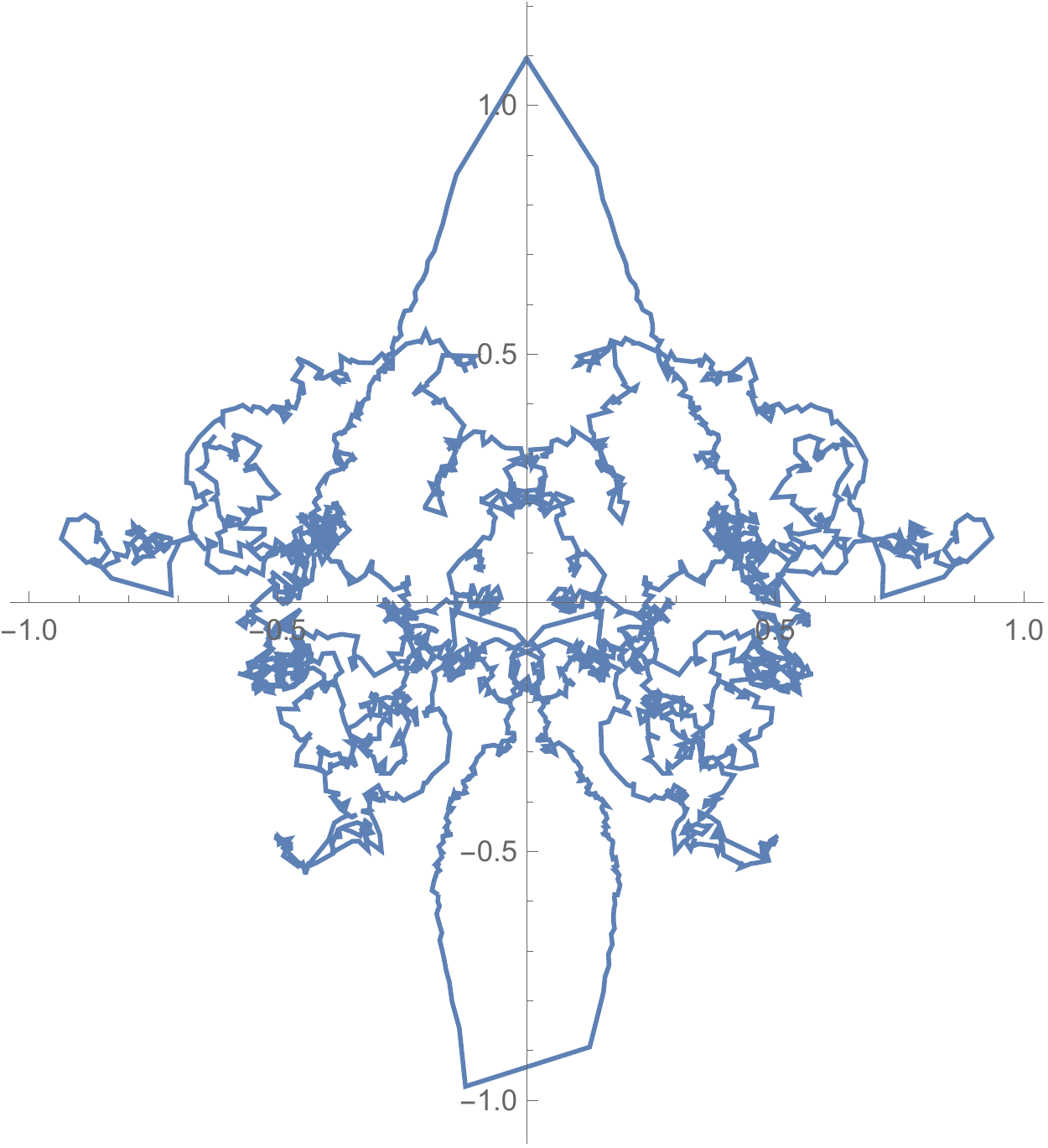}}}
\caption{Our Lissajous curves with different kind of finite Fourier sums with prime number frequencies. 
From left to right, the sum of the first 100 and 1000 second neighboring primes were taken.}	
\label{Liss3}       % Give a unique label
\end{figure}
 
%%%%%%%%%%%%%%%%%%%%%%%%%%%%%
%\section{Summary and Outlook} 
%present address: 

%%%%%%%%%%%%%%%%%%%%%%%%%%%%%%%%%%%%%%%%%%%%%%%%%%%%%%%%%%%%%%%%%%%%%                                    
%\end{multicols}
%%%%%%%%%%%%%%%%%%%%%%%%%%%%%%%%%%%%%%%%%%%%%%%%%%%%%%%%%%%%%%%%%%%%%%%                                    

\begin{references}
\bibitem{schleich} H. Maier and W.P. Schleich, {\it{Prime Numbers 101}}, Wiley Interscience 2009. 
\bibitem{hodge} T. Estermann, {\it{Introduction to Modern Prime Number Theory}}, Cambridge the University Press, 1952. 
\bibitem{laser} I.F. Barna, M.A. Pocsai and S. Varr\'o, Eur. Phys. J. Appl. Phys. {\bf{84}}, 20101 (2018).
\bibitem{flow} I.F. Barna, M.A. Pocsai and L. M\'aty\'as, Fluid. Dyn. Res. {\bf{52}}, 015515  (2020).
\bibitem{ulamsp} M. Stein  and S.M. Ulam, "An Observation on the Distribution of Primes." American Mathematical Monthly  {\bf{74}}, 43 (1967). 
\bibitem{lawr} D. Lawrence, {\it{A Cataloge of Special Plane Curves}}, (Page 178.) Dover Publication 1972.
\end{references}
\end{document}